\documentclass[11pt]{article}
\usepackage[utf8]{inputenc}
\usepackage{amsfonts}
\usepackage{graphicx}
\usepackage{amsmath}
\usepackage{bbm}
\usepackage{hyperref}
\usepackage{cleveref}
\usepackage{mathtools}
\usepackage[affil-it]{authblk}
\usepackage[dvipsnames]{xcolor}
\usepackage{caption}
\usepackage{placeins}
\usepackage{diagbox}
\usepackage{tabularx}
\usepackage{cancel}
\usepackage[symbol]{footmisc}
\usepackage{ntheorem}
\usepackage[T1]{fontenc}
\usepackage{textcomp}
\usepackage[
backend=biber,
style=numeric,
sorting=nyt
]{biblatex}

\usepackage{scalerel}
\usepackage{stackengine,wasysym}

\usepackage{geometry}
\numberwithin{equation}{section}%

\usepackage{setspace}
\hypersetup{hidelinks}

\newtheorem{theorem}{Theorem}

\newcommand\modif[1]{\color{black} #1 \color{black}}

\newcommand{\E}{\mathbbm{E}}
\newcommand{\R}{\mathbbm{R}}
\newcommand{\bcdot}{\ensuremath{\boldsymbol{\cdot}}}
\newcommand{\Av}{A^{\scriptscriptstyle v}}
\newcommand{\AT}{A^{\scriptscriptstyle T}}
\newcommand{\AS}{A^{\scriptscriptstyle S}}
\newcommand{\uS}{u_{\scriptscriptstyle S}}
\newcommand{\US}{U_{\scriptscriptstyle S}}

\newcommand{\tr}{{\scriptscriptstyle \rm T}}

\usepackage{amsmath,array}

\title{Some properties of a non-hydrostatic stochastic oceanic primitive equations model}

\author[1]{Arnaud Debussche}
\author[2]{Étienne Mémin}
\author[,2]{Antoine Moneyron \thanks{Corresponding author: antoine.moneyron@inria.fr}}

\affil[1]{Univ Rennes, CNRS, IRMAR UMR 6625, F35000, Rennes, France.}
\affil[2]{Univ Rennes, INRIA, IRMAR UMR 6625, F35000, Rennes, France.}

\date{\today}

\begin{document}

\maketitle
\thispagestyle{empty}

\begin{abstract}

\begin{center}
    \rule{0.5\linewidth}{1pt}
\end{center}
\vspace{2mm}


In this paper, we study how relaxing the classical hydrostatic balance hypothesis affects theoretical aspects of the LU primitive equations well-posedness. We focus on models that sit between incompressible 3D LU Navier-Stokes equations and standard LU primitive equations, aiming for numerical manageability while capturing non-hydrostatic phenomena. Our main result concerns the well-posedness of a specific stochastic interpretation of the LU primitive equations. This holds with rigid-lid type boundary conditions, and when the horizontal component of noise is independent of $z$, see \cite{AHHS_2022, AHHS_2023_preprint}. In fact these conditions can be related \modif{to the dynamical regime in which the primitive equations remain valid.} Moreover, under these conditions, we show that the LU primitive equations solution tends toward the one of the deterministic primitive equations for a vanishing noise, thus providing a physical coherence to the LU stochastic model. 

\begin{center}
    \rule{0.5\linewidth}{1pt}
\end{center}

\end{abstract}



\section{Introduction}

Stochastic modelling for large-scale fluid dynamics is essential but challenging due to the intrinsic complexity of geophysical flows, as they are chaotic systems with fully developed turbulence. \modif{These features induce computational limitations, which classically require to use physical approximations. However,} in the last years, stochastic modelling has emerged as a powerful setting for deriving suitable representations \cite{Berner-Al_2017, FO_2017, FOBWL_2015, MTV_1999}, allowing greater variability than deterministic large-scale representation. These models aim for plausible forecasts together with efficient uncertainty quantification.

Specifically, the location uncertainty approach (LU) has been developed in the past decade for deriving reliable stochastic models, using stochastic principles for mass, momentum, and energy conservation \cite{Mémin_2014, TML_2023}. It is applied successfully in geophysical and reduced-order models, and LU versions of fluid flow dynamics models have shown promising properties -- on classical geophysical models \cite{BCCLM_2020, RMC_2017_Pt1, RMC_2017_Pt2, RMC_2017_Pt3}, stochastic reduced order models \cite{RMHC_2017, RPMC_2021, TCM_2021} and large eddy simulation models \cite{CMH_2020, CHLM_2018, HM_2017}. The LU formalism is based on transport noises, which are intensively investigated by the mathematical community \cite{AHHS_2022, AHHS_2023_preprint, BS_2021, CFH_2019, DHM_2023, FGL_2021, FGP_10, FL_2021, GCL_2023, LCM_2023, MR_2005}.

Furthermore, the deterministic primitive equations -- which assume hydrostatic equilibrium -- are standard yet limited \cite{Vallis2017}, especially in representing important phenomena for climate such as deep convection. Remarkably, this model is known to be well-posed with rigid-lid boundary conditions \cite{CT_2007}. Moreover, a recent work showed that a stochastic representation of primitive equations with transport noise is well-posed under a ``deterministic-like'' hydrostatic hypothesis, using water world type boundary conditions \cite{AHHS_2022}. However, the authors assumed that the horizontal noise is independent of the vertical axis, which makes the barotropic and baroclinic noises tractable. In this paper, we explore \modif{weaker} hydrostatic assumptions in the LU representation of the primitive equations to model non-hydrostatic phenomena more accurately.

Our paper investigates the well-posedness of these LU primitive equations under modified hydrostatic balance, offering potential models bridging stochastic \modif{primitive equations} and 3D Navier-Stokes equations. Our main result concerns the global solutions and continuity under specific noise assumptions, aiming to better represent non-hydrostatic phenomena.

The paper is organised as follows: first we detail the assumptions made in the LU framework, then we define the function spaces to derive the abstract mathematical problem. Finally, we state our results on existence and uniqueness of solutions for the model we proposed. In particular, this model admits global pathwise solutions under regular enough conditions, and with an additional noise structure assumption. We only give a sketch of proof for our main result; a more complete one will be submitted subsequently.

\section{Oceans dynamics models in the LU framework}

The LU formulation separates the flow displacement into large-scale dynamics and highly oscillating unresolved motion, as follows,
\begin{equation}
    dX_t = u(X_t,t) dt + \sigma(X_t,t) dW_t.
\end{equation}
Here, $X$ denotes the \modif{ (3D) \emph{Lagrangian}} displacement, defined in a bounded cylindrical domain $\mathcal{S} = \mathcal{S}_H \times [-h,0] \subset \R^3$, where $\mathcal{S}_H$ is a subset of $\R^2$ with smooth boundary. \modif{In this formalism, $u : \mathcal{S} \times [0,T] \rightarrow \R^3$ is the \emph{Eulerian} velocity of the fluid flow.} The large-scale velocity $u(X_t,t)$ correlates in both space and time, while the unresolved small-scale velocity $\sigma(X_t,t) dW$ is uncorrelated-in-time yet correlated-in-space. We further refer to this second component as a noise term, which \modif{must be} interpreted in the It\={o} sense.

\modif{Let us define this noise term more precisely: consider a cylindrical Wiener process $W$ on the space of square integrable functions $\mathcal{W} := L^2(\mathcal{S},\R^3)$. Thus, there exists a Hilbert orthonormal basis $(e_i)_{i \in \mathbbm{N}}$ of $\mathcal{W}$ and a sequence of independent standard Brownian motions $(\hat{\beta}^i)_{i \in \mathbbm{N}}$ on a filtered probability space $(\Omega, \mathcal{F}, (\mathcal{F}_t)_t, \mathbbm{P})$ such that,
\begin{equation*}
    W=\sum_{i \in \mathbbm{N}} \hat{\beta}^i e_i.
\end{equation*}
Note that the sum $\sum_{i \in \mathbbm{N}} \hat{\beta}^i e_i$ does \emph{not} converge \emph{a priori} in $\mathcal{W}$. Hence we interpret the previous identity in a space $\mathcal{U}$ including $\mathcal{W}$, such that the embedding $\mathcal{W} \hookrightarrow \mathcal{U}$ is Hilbert-Schmidt. Typically, $\mathcal{U}$ can be the dual space of any reproducing kernel Hilbert subspace of $\mathcal{W}$ for the inner product $(\bcdot, \bcdot)_{\mathcal{W}}$, that is for instance $H^{-s}(\mathcal{S}, \R^3)$ with $s > \frac{3}{2}$.

Then, we define the noise through a deterministic time-dependent correlation operator $\sigma_t$: let $\hat{\sigma} : [0,T] \rightarrow L^2(\mathcal{S}^2,\R^3)$ a bounded symmetric kernel, and define
\begin{equation*}
    (\sigma_t f)(x) = \int_{\mathcal{S}} \hat{\sigma}(x,y,t) f(y) dy, \quad \forall f \in \mathcal{W}.
\end{equation*}
With this definition, $\sigma_t$ is a Hilbert-Schmidt operator mapping $\mathcal{W}$ into itself, so that the noise can be defined as
\begin{equation*}
    \sigma_t W_t = \sum_{i \in \mathbbm{N}} \hat{\beta}_t^i \sigma_t e_i,
\end{equation*}
where the previous series converges in the sense of $L^2(\Omega, \mathcal{W})$. Here we interpret $\mathcal{W}$ as the space carrying the process $\sigma_t W_t$, while the notation $L^2(\mathcal{S},\R^3)$ is kept for denoting the space of tridimensional velocities.  Importantly, $\sigma_t W_t$ is an abuse of notation as $\sigma_t$ is an operator is defined on $\mathcal{W}$ while $W_t$ does not converge on $\mathcal{W}$ -- that is to say ``$\sigma_t (W_t)$'' \emph{is not} properly defined.

Moreover, we consider a family of eigenfunctions $(\phi_k)_k$ of the operator $\sigma_t$, which we scale by their corresponding eigenvalues, and such that $(\phi_k)_k$ is a Hilbert basis of $\mathcal{W}$. By writing the previous series in terms of $(\phi_k)_k$, it can be shown that there exists a sequence of independent standard Brownian motions $(\beta_t^k)_k$, defined on the previously introduced filtered space $(\Omega, \mathcal{F}, (\mathcal{F}_t)_t, \mathbbm{P})$, so that
\begin{equation*}
    \sigma_t W_t = \sum_{k \in \mathbbm{N}} \beta_t^k \phi_k(t).
\end{equation*}
As such, $(\Omega, \mathcal{F}, (\mathcal{F}_t)_t, \mathbbm{P}, W)$ is a stochastic basis. In addition, the previous series converges in $\mathcal{W}$ almost surely, and in the sense of $L^p(\Omega, \mathcal{W})$ for all $p \in \mathbbm{N}$ \cite{book_DPZ_2014, book_GM_2010}.}


\bigskip

\modif{Furthermore, we may define} the variance tensor as follows, which is the diagonal part of the covariance tensor,
\begin{equation}
    a(x,t) = \int_{\mathcal{S}} \hat{\sigma}(x,y,t) \hat{\sigma}(y,x,t) dy = \sum_{k=0}^\infty \phi_k(x,t) \phi_k(x,t)^\tr.
\end{equation}
Let us note that, in full generality \modif{-- that is when $\hat{\sigma}_t$ is itself a random function --} the operator-valued process $\sigma_t$ is subject to an integrability condition 
\modif{$$\mathbbm{P}\Big(\int_0^T \|a(\bcdot,t)\|_{L^2(\mathcal{S}, \R^{ 3\times 3})} dt <\infty\Big) =1,$$
where $\| \bcdot \|_{L^2(\mathcal{S}, \R^{ 3\times 3}) }$ is the Hilbert norm associated to $L^2(\mathcal{S}, \R^{ 3\times 3})$, the matrix space $ \R^{ 3\times 3}$ being equipped with the Frobenius norm. As such, the integral $\int_0^t \sigma_s dW_s $ is a $\mathcal{W}$-valued Gaussian process with expectation zero and bounded variance: $\E\big[\|\int_0^t \sigma_s dW_s\|_{L^2}^2\big]<\infty$.} The quadratic variation of $\int_0^t (\sigma_t d W_s)(x)$ is given by the finite variation process $\int_0^t a(x,s) ds.$

In a similar way as deriving the classical Navier-Stokes equations, the LU Navier-Stokes equations emerge through a stochastic version of the Reynolds transport theorem (SRTT) \cite{Mémin_2014}. Let $q$ be a random scalar, within a volume $\mathcal{V}(t)$ transported by the flow. Then, for incompressible unresolved flows -- that is $\nabla \bcdot \sigma_t=0$ -- the SRTT \modif{reads}
\begin{gather}
    d\Big( \int_{\mathcal{V}(t)} q(x,t) dx \Big) = \int_{\mathcal{V}(t)} \big(\mathbbm{D}_t q + q \nabla \bcdot (u-\uS) dt\big) dx,\\
    \mathbbm{D}_t q = d_t q + (u - \uS) \bcdot \nabla q \: dt + \sigma dW_t \bcdot \nabla q - \frac{1}{2} \nabla \bcdot (a \nabla q) dt, \label{transport-operator}
\end{gather}
involving an additional drift $\uS = \frac{1}{2} \nabla \bcdot a$, coined as the It\={o}-Stokes drift in \cite{BCCLM_2020}. \modif{Consequently, $\uS$ is a vector field from $\mathcal{S} \times [0,T]$ to $\R^3$, alike any trajectory of $u$. In addition,} $d_t q(x,t) = q(x,t+dt) - q(x,t)$ is the forward time increment at a fixed spatial point $x$, and $\mathbbm{D}_t q$ is a stochastic transport operator introduced in \cite{Mémin_2014, RMC_2017_Pt1}, playing the role of the material derivative. The It\={o}-Stokes drift is directly related to the divergence of the variance tensor $a$, representing the effects of noise inhomogeneity on large-scale dynamics. Such advection terms are often added in \modif{large-scale ocean dynamics models} -- under the name of bollus velocity -- to take account for surface waves and Langmuir turbulence \cite{CL_1976, GOP_2004, MSM_1997}. As shown in \cite{BCCLM_2020}, the LU framework holds similar features which accounts for the effects of the small-scale inhomogeneity on the large-scale flow. Additionally, the stochastic transport operator represents physically interpretable terms for large-scale flows representation. Namely, \modif{equation \eqref{transport-operator} gathers the following four terms}
\begin{itemize}
    \item an evolution term $d_t q$,
    \item a transport term $(u - \uS) \bcdot \nabla q \: dt$, that is the advection of the large-scale quantity $q$ by the large-scale It\={o}-Stokes drift corrected velocity $u-\uS$,
    \item a transport term $\sigma dW_t \bcdot \nabla q$, that is the advection of $q$ by the unresolved velocity,
    \item an inhomogeneous diffusion term $- \frac{1}{2} \nabla \bcdot (a \nabla q) dt$ representing small-scale mixing effects on $q$.
\end{itemize}
Remarkably, the energy associated to the backscattering term $\sigma dW_t \bcdot \nabla q$ is exactly compensated by the stochastic diffusion term $- \frac{1}{2} \nabla \bcdot (a \nabla q) dt$ \cite{RMC_2017_Pt1}. This equilibrium can be interpreted as an instance of a fluctuation-dissipation theorem.

\subsection{The LU primitive equations}

Let us derive the LU primitive equations. Assuming that the flow is isochoric with constant material density, and that the noise is divergence-free alongside with a divergence-free corresponding It\={o}-Stokes drift, i.e.
\begin{equation}
    \nabla \bcdot u = \nabla \bcdot \uS =0, \quad \nabla \bcdot \sigma_t dW_t =0,
\end{equation}
we deduce that $\mathbbm{D}_t q =0$, for any conservative scalar quantity $q$. In addition, we model density through a linear law of state involving salinity and temperature,
\begin{equation}
    \rho = \rho_0\Big(1 + \beta_T (T-T_r) + \beta_S (S-S_r)\Big),
\end{equation}
with $\rho_0$ the reference density of the ocean at a typical temperature $T_r$ and salinity $S_r$, assuming that the thermodynamic parameters $\beta_T := \frac{1}{\rho_0} \frac{\partial \rho}{\partial T}$ and $\beta_S := \frac{1}{\rho_0} \frac{\partial \rho}{\partial T}$ are constant. Plus, we define three anisotropic diffusion operators that is used further: for $i \in \{v,T,S\}$, provided that the viscosities $\mu_i, \nu_i$ are given \emph{a priori}, denote
\begin{align}
    A^{\scriptscriptstyle i} &= -\mu_i (\partial_{xx} + \partial_{yy}) - \nu_i  \partial_{zz}.
\end{align}

Applying the SRTT to the conservation of momentum principle in rotating frame -- see \cite{DHM_2023} in non rotating frame -- we derive the following stochastic equations of motion,
\begin{equation}
    \mathbbm{D}_t u + f k \times (u \: dt + \sigma dW_t) = -\frac{1}{\rho_0} \nabla (p \: dt + dp_t^\sigma) - \Av (u \: dt + \sigma dW_t). \label{NS-u}
\end{equation}
Importantly, we have introduced a martingale noise pressure term $dp_t^\sigma$ in addition to the classical pressure term $p \: dt$, due to stochastic modelling. Applying the SRTT similarly to the conservation of energy and saline mass, we find
\begin{align}
    \mathbbm{D}_t T = - \AT T, \quad \mathbbm{D}_t S = - \AS S. \label{NS-TS}
\end{align}

Now, define $u^* = u - \uS$, and note that by assumption $u^*$, $u$ and $\uS$ are divergence free. Therefore, if we denote by $v^*$, $v$ and $v_s$ their respective horizontal components, we can express their vertical components via the integro-differential operator
\begin{equation}
    w(\bcdot) = \int_z^0 \nabla_H \bcdot (\bcdot), \label{w-definition}
\end{equation}
under the hypothesis that $w=0$ when $z = 0$ or $z=-h$. The horizontal gradient operator is denoted by $\nabla_H = (\partial_x \: \partial_y)^\tr$, and we use further $\Delta_H(\bcdot) = \nabla_H \bcdot (\nabla_H (\bcdot) )$ to denote the horizontal Laplace operator. Also, write $\sigma^H dW_t$ and $\sigma^w dW_t$ the horizontal and vertical components of $\sigma dW_t$, respectively. Thus, the horizontal and vertical momentum equations \modif{read}
\begin{gather}
    \mathbbm{D}_t v + \Gamma (v \: dt + \sigma^H dW_t) =  -\Av (v\:dt  + \sigma^H dW_t) - \frac{1}{\rho_0} \nabla_H (p \: dt + dp_t^\sigma), \label{NS-v}\\
    \mathbbm{D}_t w = -\Av (w\:dt  + \sigma^w dW_t) - \frac{1}{\rho_0} \partial_z (p \: dt + dp_t^\sigma) - \frac{\rho}{\rho_0} g dt, \label{NS-w}
\end{gather}
where $\Gamma((a \quad b)^\tr) = f (-b \quad a)^\tr$ is the horizontal projection of the Coriolis term.

A recent study -- see article \cite{AHHS_2022} -- derived stochastic primitive equations by considering a ``deterministic-like'' hydrostatic hypothesis: assuming that the vertical acceleration is negligible compared to the gravity term, the vertical momentum equation boils down to
\begin{equation}
    \partial_z p + \rho g = 0, \quad \text{and} \quad \partial_z dp_t^\sigma = 0. \label{strong-hydro}
\end{equation}
We further call this assumption the \emph{strong hydrostatic hypothesis}. The validity of this hydrostatic balance corresponds to a regime with  small ratio between the squared aspect ratio $\alpha^2=h^2/L^2$ -- with $h$ and $L$ denoting vertical and horizontal scales, respectively -- and the Richardson number $Ri = \frac{N^2}{(\partial_z v)^2} \sim N^2 h^2/U^2$, defined from the ratio between the stratification given by the Brunt-Väsäilä frequency $N^2= - \frac{g}{\rho_0} \partial_z \rho $ and the squared vertical shear of the horizontal velocity \cite{MHPA_1997}. In the stochastic setting, the strong hydrostatic balance holds if the noise is small enough not to unbalance this regime. 

Gathering all the equations described previously and assuming that the strong hydrostatic hypothesis holds, we obtain the following problem
\begin{gather*}
    \mathbbm{D}_t v + \Gamma (v \: dt + \sigma^H dW_t) =  -\Av (v\:dt  + \sigma^H dW_t) - \frac{1}{\rho_0} \nabla_H (p \: dt + dp_t^\sigma),\\
    \mathbbm{D}_t T = - \AT T dt, \quad
    \mathbbm{D}_t S = - \AS S dt,\\
    \nabla_H \bcdot v + \partial_z w =0,\\
    \partial_z p + \rho g = 0, \quad 
    \partial_z dp_t^\sigma = 0,\\
    \rho = \rho_0(1 + \beta_T (T-T_r) + \beta_S (S-S_r)).    
\end{gather*}
This system is closely aligned with the deterministic primitive system, since the stochastic transport operator is interpreted as a material derivative. The aforementioned model studied in \cite{AHHS_2022} corresponds to the model above, with the stochastic diffusion term $\frac{1}{2} \nabla \bcdot (a \nabla(\bcdot))$ being replaced by $\nu_\sigma \Delta (\bcdot)$, where $\nu_\sigma>0$ is a constant.Well-posedness results were established assuming that the initial condition is smooth enough, with periodic horizontal boundary conditions and rigid-lid type vertical boundary conditions -- termed \emph{water world} or \emph{aqua planet}. Their key results highlighted local-in-time well-posedness, and global-in-time well-posedness when the horizontal component of the noise is barotropic, i.e. independent of the z-coordinate. We are aiming to extend these well-posedness proofs to more general models, introducing different assumptions about the vertical momentum equation, which we term \emph{weak hydrostatic hypotheses}. These assumptions accommodate relaxed hydrostatic equilibria, touching regimes at the edge of the deterministic hydrostatic assumption validity. This is expected to capture better non-hydrostatic phenomena (like wind or buoyancy-driven turbulence), and deep oceanic convection, influenced by strong noise disrupting the strong hydrostatic regime. We further present one possible interpretation of the \emph{weak hydrostatic hypothesis}.

To derive our new model, we neglect only the large-scale contribution of the vertical acceleration in the vertical momentum equation \eqref{NS-w},
\begin{align*}
    - &\frac{1}{\rho_0} \partial_z (p \: dt + dp_t^\sigma) - \frac{\rho}{\rho_0} g dt = \mathbbm{D}_t w + \Av(w\:dt + \sigma^w dW_t)\\ &= \underbrace{d_t w + u \bcdot \nabla w \: dt + \Av(w\:dt)}_{\approx 0} - \uS \bcdot \nabla w dt + \sigma dW_t \bcdot \nabla w - \frac{1}{2} \nabla \bcdot (a \nabla w) dt + \Av(\sigma^w dW_t)\\
    &\approx - \uS \bcdot \nabla w dt + \sigma dW_t \bcdot \nabla w - \frac{1}{2} \nabla \bcdot (a \nabla w) dt + \Av(\sigma dW_t).
\end{align*}
that is,
\begin{equation}
    \frac{1}{\rho_0} \partial_z p = - \frac{\rho}{\rho_0} g + \frac{1}{2} \nabla \bcdot (a \nabla w) + \uS \bcdot \nabla w dt, \text{ and } \quad
    \frac{1}{\rho_0} \partial_z dp_t^\sigma = -\sigma dW_t \bcdot \nabla w - \Av(\sigma^w dW_t). \label{weak-NS-hydro}
\end{equation}
Notice that the stochastic advection and diffusion terms have been kept in this interpretation of the hydrostatic hypothesis. \modif{In its deterministic version,} this hypothesis neglects the whole vertical acceleration term and molecular diffusion term, arguing that they are negligible compared to the gravity term $\frac{\rho}{\rho_0} g$. Transitioning to the stochastic case, we claim the same for the large-scale vertical terms $d_t w + u \bcdot \nabla w dt$ and $A^v w dt$. However, due to stochastic modelling, three terms remain: the back scattering noise advection term $\sigma dW_t \bcdot \nabla w$, the It\={o}-Stokes drift \modif{advection} term $\uS \bcdot \nabla w dt$, and the diffusion term $\frac{1}{2} \nabla \bcdot (a \nabla w) dt$. Differently from $u \bcdot \nabla w dt$, \modif{the term} $\sigma dW_t \bcdot \nabla w$ cannot be neglected compared to $\frac{\rho}{\rho_0} g dt$ since the latter is a pure bounded variation term. In addition, the terms $\uS \bcdot \nabla w dt$ and $\frac{1}{2} \nabla \bcdot (a \nabla w) dt$ are kept since they depend directly on the noise -- if $\Upsilon^{1/2}$ is the noise scaling, then they scale like $\Upsilon$. This formulation leads to the following problem
\begin{gather*}
    \mathbbm{D}_t v + \Gamma (v \: dt + \sigma^H dW_t) =  -\Av (v\:dt  + \sigma^H dW_t) - \frac{1}{\rho_0} \nabla_H (p \: dt + dp_t^\sigma),\\
    \mathbbm{D}_t T = - \AT T dt, \quad
    \mathbbm{D}_t S = - \AS S dt,\\
    \nabla_H \bcdot v + \partial_z w =0,\\
    \frac{1}{\rho_0} \partial_z p + \frac{\rho}{\rho_0} g = \frac{1}{2} \nabla \bcdot (a \nabla w) +  \uS \bcdot \nabla w, \quad
    \frac{1}{\rho_0} \partial_z dp_t^\sigma = -\sigma dW_t \bcdot \nabla w - \Av(\sigma^w dW_t),\\
    \rho = \rho_0(1 + \beta_T (T-T_r) + \beta_S (S-S_r)).    
\end{gather*}
Considering such a problem, a major difficulty arises: the presence of the transport noise $\sigma dW_t \bcdot \nabla w$ in the stochastic pressure term induces the following terms in horizontal velocity dynamics,
$$\nabla_H \int_z^0 \sigma dW_t \bcdot \nabla w \: dz', \quad \nabla_H \int_z^0 \uS \bcdot \nabla w \: dz' \: dt.$$
We don't know how to make a suitable energy estimate for global pathwise existence with such terms, notably because three derivatives on the horizontal velocity $v$ are involved -- since $w(v) = \int_z^0 \nabla_H \bcdot v$. To tackle this issue, we propose one way to regularise the previous system.

\paragraph{Weak low-pass filtered hydrostatic hypothesis:}

To enforce greater regularity for the vertical transport noise term, we define a regularising convolution kernel $K$, and replace equation \eqref{weak-NS-hydro} by
\begin{equation}
    \frac{1}{\rho_0} \partial_z p = - \frac{\rho}{\rho_0} g + \frac{1}{2} \nabla \bcdot (a^K \nabla w) + K*[\uS \bcdot \nabla w], \quad
    \frac{1}{\rho_0} \partial_z dp_t^\sigma = -K*[\sigma dW_t \bcdot \nabla w] - \Av(\sigma^w dW_t). \label{weak-split}
\end{equation}
In the previous formulae, we denote by $a^K$ the operator $f \longmapsto \sum_{k=0}^\infty \phi_k \mathcal{C}_K \mathcal{C}_K^* (\phi_k^\tr f)$, and by $\mathcal{C}_K$ the operator $f \longmapsto K*f$. The regularising kernel only affects the vertical transport noise, and not potential vertical additive noises. Also, the stochastic diffusion operator $\frac{1}{2} \nabla \bcdot (a^K \nabla (\bcdot))$ is chosen to be the covariation correction term associated to $K*[\sigma dW_t \bcdot \nabla w]$. We refer to this assumption as the \emph{weak low-pass filtered hydrostatic hypothesis}.

This approach consists in filtering the vertical transport noise and disregarding the vertical acceleration of the resolved component of velocity. The noise terms, alongside the stochastic diffusion term, account for the deviation from a strong hydrostatic equilibrium. By convolving the vertical transport noise with $K$, we effectively truncate its highest frequencies. This new hypothesis relaxes the strong hydrostatic balance, allowing for more general stochastic pressures and extending the validity of the system dynamical regime beyond the strong hydrostatic case. Comparing this methodology to the one presented in \cite{AHHS_2023_preprint}, authors introduced a temperature noise impacting the pressure equation, that is a perturbation of \emph{thermodynamic} origin. Unlike theirs, our model involves transport noise of the vertical velocity component, that is a perturbation of \emph{mechanical} origin. This approach retains more terms linked to the vertical velocity $w$, accounting for the influence of unresolved small-scale velocity -- like turbulence or submesoscale components -- on the large-scale vertical velocity. Hence, the structure of this problem resembles a more genuine tridimensional problem, justifying the additional regularization through a filtering kernel. Using a filtering kernel is a common practice in defining numerical models for primitive equations. It is also prevalent in establishing well-posedness for specific subgrid models like the Gent-McWilliam model, particularly in the field of mesoscale dynamics, as highlighted in \cite{KT_2023}. Thus, assuming that this weak hydrostatic hypothesis holds rather than the strong one, one derives the following problem
\begin{gather}
    \mathbbm{D}_t v + \Gamma (v \: dt + \sigma^H dW_t) =  -\Av (v\:dt  + \sigma^H dW_t) - \frac{1}{\rho_0} \nabla_H (p \: dt + dp_t^\sigma), \label{primitive-v}\\
    \mathbbm{D}_t T = - \AT T dt,\label{primitive-T}\\
    \mathbbm{D}_t S = - \AS S dt,\label{primitive-S}\\
    \nabla_H \bcdot v + \partial_z w =0,\\
    \frac{1}{\rho_0}\partial_z p + \frac{\rho}{\rho_0} g = \frac{1}{2} \nabla \bcdot (a^K \nabla w) + K*[\uS \bcdot \nabla w],\label{primitive-p}\\
    \frac{1}{\rho_0}\partial_z dp_t^\sigma +K*[\sigma dW_t \bcdot \nabla w ]+ \Av(\sigma^w dW_t) = 0,\label{primitive-dp}\\
    \rho = \rho_0\Big(1 + \beta_T (T-T_r) + \beta_S (S-S_r)\Big).    
\end{gather}
In addition, decompose the boundary as $\partial S = \Gamma_u \cup \Gamma_b \cup \Gamma_l$  -- respectively the upper, bottom and lateral boundaries -- and equip this problem with the following rigid-lid type boundary conditions \cite{BS_2021,DGHT_2011}
\begin{align}
    \partial_z v  = 0, && w = 0, && \nu_T\partial_z T + \alpha_T T =0, && \partial_z S=0 && \text{on } \Gamma_u, \label{boundary-conditions}\\
    v = 0, && w = 0, && \partial_z T =0, && \partial_z S=0 && \text{on } \Gamma_b, \nonumber \\
    v = 0, &&  && \partial_{\vec{n}_H} T =0, && \partial_{\vec{n}_H} S =0 && \text{on } \Gamma_l, \nonumber
\end{align}
and initial conditions
\begin{align}
    v(t=0) = v_0 \in H^1(\mathcal{S},\R^2), && T(t=0) = T_0 \in H^1(\mathcal{S},\R), &&
    S(t=0) = S_0 \in H^1(\mathcal{S},\R),
    \label{initial-conditions}
\end{align}
fulfilling the previous boundary conditions. We further assume that the noise, its variance tensor and their gradients cancel on the horizontal boundary $\Gamma_l$. Also, we assume that the vertical noise $\sigma_z dW_t$ cancels on the vertical upper and bottom boundaries $\Gamma_u \cup \Gamma_b$. Furthermore, the noise is assumed to be regular enough in the following sense,
\begin{align}
    \sup_{t\in [0,T]} \sum_{k=0}^\infty \|\phi_k\|_{H^4(\mathcal{S}, \R^3)}^2 < \infty, \label{smoothness-noise} && \uS \in L^\infty\Big([0,T], H^4(\mathcal{S}, \R^3)\Big), \\
    d_t \uS \in L^\infty\Big([0,T],H^3(\mathcal{S}, \R^3)\Big), && a \nabla \uS \in L^\infty\Big([0,T],H^2(\mathcal{S},\R^3)\Big). \nonumber
\end{align}
and 
\begin{equation}
    a \in H^1\Big([0,T],H^3(\mathcal{S},\R^{3 \times 3})\Big). \label{smoothness-noise-2}
\end{equation}
These regularity assumptions are not limiting in practice \modif{because} most ocean models consider spatially smooth noises \cite{TML_2023}, since these are the physically observed ones. Also, remark that the weak low-pass filtered hydrostatic hypothesis adds a stochastic contribution to the pressure terms compared to the strong hypothesis -- see equations \eqref{primitive-p} and \eqref{primitive-dp} --
\begin{multline*}
    \frac{1}{\rho_0}\partial_z p + \frac{\rho}{\rho_0} g = \frac{1}{2} \nabla \bcdot (a^K \nabla w) + K*[\uS \bcdot \nabla w], \\
    \text{and } \quad \frac{1}{\rho_0}\partial_z dp_t^\sigma + K * [\sigma dW_t \bcdot \nabla w] + \Av(\sigma^w dW_t) = 0.
\end{multline*}
Therefore, this impacts the horizontal momentum equation via the horizontal pressure gradients, so one has
\begin{align*}
    \frac{1}{\rho_0}\nabla_H p &= - \nabla_H \int_z^0 \frac{\rho}{\rho_0} g \: dz'  +  \nabla_H \int_z^0 \frac{1}{2} \nabla \bcdot (a^K \nabla w) + K*[ \uS \bcdot \nabla w ] \: dz' + \frac{1}{\rho_0} \nabla_H p^s\\
    &=: - \nabla_H \int_z^0 \frac{\rho}{\rho_0} g \: dz' + \frac{1}{\rho_0}\nabla_H p^{weak} + \frac{1}{\rho_0} \nabla_H p^s,\\
    \frac{1}{\rho_0}\nabla_H dp_t^\sigma &=  -\nabla_H \int_z^0 K*[\sigma dW_t \bcdot \nabla w] - \Av(\sigma^w dW_t) \: dz' + \frac{1}{\rho_0} \nabla_H dp^{\sigma,s}_t\\
    &=: \frac{1}{\rho_0}\nabla_H dp^{\sigma, weak}_t + \frac{1}{\rho_0} \nabla_H dp^{\sigma,s}_t,
\end{align*}
where $p^{weak}$ and $p^{\sigma, weak}$ denote respectively the bounded variation and martingale pressures associated to the weakening of the hydrostatic hypothesis. Introducing new additive noise terms does not pose a challenge for a sufficiently regular $\sigma dW_t$. However, the key difficulty arises from the presence of the term $\nabla_H \int_z^0 \sigma dW_t \bcdot \nabla w \: dz'$ in the expression of $\frac{1}{\rho_0} \nabla_H dp^{\sigma}_t$. This term represents the horizontal impact of the vertical transport noise $\sigma dW_t \bcdot \nabla w$. Regularising this term using a smoothing filter inherently enlarges its spatial scale, causing the spatial scale of the vertical transport noise to remain above the resolution cutoff scale without aliasing artifacts. Likewise, $\nabla_H \int_z^0 \nabla \bcdot (a \nabla w) \: dz'$ represents the horizontal influence of the covariation correction originating from the LU Navier-Stokes equations. Establishing the well-posedness of such a model, subject to appropriate regularisation and structural conditions, stands as our primary outcome.

\subsection{Definition of the spaces}\label{subsec-def-spaces}

In this subsection, we define the function spaces that are used further. Remind that the spatial domain is denoted by $\mathcal{S} = \mathcal{S}_H \times [-h,0] \subset \R^3$. First, we introduce the following inner products
\begin{gather*}
    (v,v^\sharp)_{H_1} = (v,v^\sharp)_{L^2(\mathcal{S},\R^2)}, \quad (v,v^\sharp)_{V_1} = (\nabla v, \nabla v^\sharp)_{L^2(\mathcal{S},\R^2)}\\
    (T,T^\sharp)_{H_2} = (T,T^\sharp)_{L^2(\mathcal{S},\R)}, \quad (T,T^\sharp)_{V_2} = (\nabla T,\nabla T^\sharp)_{L^2} + \frac{\alpha_T}{\nu_T}(T,T^\sharp)_{L^2(\Gamma_u, \R)},\\
    (S,S^\sharp)_{H_3} = (S,S^\sharp)_{L^2(\mathcal{S},\R)}, \quad (S,S^\sharp)_{V_3} = (\nabla S,\nabla S^\sharp)_{L^2(\mathcal{S},\R)},
\end{gather*}
and let
\begin{align}
    (U,U^\sharp)_H = (v,v^\sharp)_{H_1} + (T,T^\sharp)_{H_2} + (S,S^\sharp)_{H_3}, \quad (U,U^\sharp)_V = (v,v^\sharp)_{V_1} + (T,T^\sharp)_{V_2} + (S,S^\sharp)_{V_3},
\end{align}
for all $U,U^\sharp \in L^2(\mathcal{S},\R^4)$. Also we denote by \modif{$\|\bcdot\|_H, \|\bcdot\|_{H_i}$ and $\|\bcdot\|_V,\|\bcdot\|_{V_i}$ the associated norms. With a slight abuse of notation, we may write $\|\bcdot\|_H, \|\bcdot\|_V$ in place of $\|\bcdot\|_{H_i}, \|\bcdot\|_{V_i}$ respectively, and similarly use $(\bcdot, \bcdot)_H, (\bcdot, \bcdot)_V$ rather than $(\bcdot, \bcdot)_{H_i}, (\bcdot, \bcdot)_{V_i}$. Moreover, notice} that we distinguished the inner products $(\bcdot, \bcdot)_{H_2}$ and $(\bcdot, \bcdot)_{H_3}$, even if they denote the same operation. This is for consistency with the following definitions of the spaces $H_2$ and $H_3$.

Then, denote by $\mathcal{V}_1$ the space of functions $C^\infty(\mathcal{S},\R^2)$ with a compact support strictly included in $\mathcal{S}$, such that for all $v \in \mathcal{V}_1$, $\nabla_H \bcdot \int_{-h}^0 v = 0$. Plus, define $\mathcal{V}_2 = \mathcal{V}_3$ the space of functions $C^\infty(\mathcal{S},\R)$. Denote by $H_i$ the closure of $\mathcal{V}_i$ for the norm $\|.\|_{H_i}$, and $V_i$ its closure by $\|.\|_{V_i}$. Eventually, define $H=H_1 \times H_2 \times H_3$ and $V=V_1 \times V_2 \times V_3$, which are also the closures of $\mathcal{V}_1 \times \mathcal{V}_2 \times \mathcal{V}_3$ by $\|.\|_{H}$ and $\|.\|_{V}$, respectively. Often, by abuse of notation, we write $(\bcdot,\bcdot)_H$ instead of $(\bcdot,\bcdot)_{V' \times V}$. More generally, if $K$ is a subspace of $H$ and $K'$ its dual space, we write $(\bcdot,\bcdot)_H$ instead of $(\bcdot,\bcdot)_{K' \times K}$. Moreover, we define $\mathcal{D}(A) = V \cap H^2(\mathcal{S},\R^4)$, where $A=(\Av \AT \AS)^\tr$ gathers the three diffusion operators. As such, $A : \mathcal{D}(A) \rightarrow H$ is an unbounded operator. In addition, for any Hilbert spaces $\mathcal{H}_1$ and $\mathcal{H}_2$, we define $\mathcal{L}_2( \mathcal{H}_1,\mathcal{H}_2)$ the space of Hilbert-Schmidt operators from $\mathcal{H}_1$ to $\mathcal{H}_2$, and $\|\bcdot\|_{\mathcal{L}_2( \mathcal{H}_1,\mathcal{H}_2)}$ its associated norm.

Again, we distinguished the spaces $H_2$ and $H_3$, even though they formally denote the same space. However, $V_2$ and $V_3$ \emph{are} different spaces since they are not equipped with the same inner products due to different boundary conditions on the temperature and salinity \modif{-- Robin and Neumann respectively.} This distinction allows to interpret $H_2$ and $V_2$ as temperature spaces, and $H_3$ and $V_3$ as salinity spaces. In addition, $H_1$ and $V_1$ are interpreted as horizontal velocity spaces \modif{-- since their elements are $\R^2$-valued processes.} Using this formalism, the vertical velocity $w$ is written as a functional of the horizontal velocity $v \in H_1$ through the continuity equation, namely $w(v) = \int_z^0 \nabla_H \cdot v \: dz'$.

Eventually, we define the barotropic and baroclinic projectors $\mathcal{A}_2 : \R^3 \rightarrow \R^2$, $\mathcal{A} : \R^3 \rightarrow \R^3$ and $\mathcal{R} : \R^3 \rightarrow \R^3$ of the velocity component as follows. For $v \in \mathcal{V}_1$, $h$ being the depth of the ocean, let
\begin{equation}
    \mathcal{A}_2 [v](x,y) = \frac{1}{h} \int_{-h}^0 v(x,y,z') dz', \quad \mathcal{A} [v](x,y,z) = \mathcal{A}_2 [v](x,y), \quad \mathcal{R} [v] = v - \mathcal{A} [v]. \label{def-barotropic-baroclinic}
\end{equation}
Remark that $\mathcal{A}$ and $\mathcal{R}$ are orthogonal projectors with respect to the inner product $(\bcdot, \bcdot )_{H}$. To simplify notations, we may use $\bar{v}$ in place of $\mathcal{A}_2 [v]$ or $\mathcal{A} [v]$, and $\tilde{v}$ in place of $\mathcal{R} [v]$.

\subsection{Abstract formulation of the problems}

In this section, we aim to express, in abstract form, the previous problem under the weak low-pass filtered hydrostatic hypothesis \eqref{weak-NS-hydro}. First, define the 4D vector $U$, representing the state of the system, and the correction $U^*$ of $U$ by the It\={o}-Stokes drift, as
$$U = (v, T, S)^\tr, \quad U^* = (v^*, T, S)^\tr.$$
Then, denote the advection operator by
\begin{equation*}
    B(U^*,\bcdot) = B(v^*,\bcdot) := (v^* \bcdot \nabla_H)( \bcdot) + w(v^*) \partial_z (\bcdot).
\end{equation*}
Moreover, define two Leray type projectors $\mathbf{P}^v$ and $\mathbf{P}$ as follows \cite{BS_2021}:
\begin{equation}
    \mathbf{P}^v(v) = \mathbf{P}_{2D} \mathcal{A} (v) + \mathcal{R} (v), \quad \mathbf{P}(U) = (\mathbf{P}^v(v), T, S)^\tr,
\end{equation}
where $\mathbf{P}_{2D}$ is the standard 2D Leray projector, which is associated to the barotropic component $\mathcal{A} (v)$. Notice that the baroclinic component $\mathcal{R} (v)$ is left unchanged by the projector $\mathbf{P}^v$, that is $\mathbf{P}$ only affects the barotropic component of velocity. For notational convenience we keep the same notations for the composition of the following operators with the Leray projector:
\modif{\begin{gather}
    AU = \mathbf{P} ( \Av v, \: \AT T, \: \AS S)^\tr,\quad CU = \mathbf{P} ( Cv, \: 0, \: 0)^\tr, \nonumber \\
    B(U^*,U) = \mathbf{P} ( B(v^*,v), \: B(v^*,T), \: B(v^*,S))^\tr. \label{def-operators}
\end{gather}}
Thus, if $U \in \mathcal{D}(A)$, then $AU \in H$, and we obtain the relation
\begin{multline}
    d_t U + \Bigg[AU + B(U^*,U) + CU + \frac{1}{\rho_0} \mathbf{P} \begin{pmatrix}
        \nabla_H p \\0 \\0
    \end{pmatrix} - \frac{1}{2} \mathbf{P} \nabla \bcdot(a\nabla U)\Bigg] dt = - \mathbf{P}( \sigma dW_t \bcdot \nabla U)\\
    - (A + \Gamma) (\sigma^H dW_t) - \frac{1}{\rho_0} \mathbf{P} \begin{pmatrix}
        \nabla_H dp_t^\sigma \\0 \\0
    \end{pmatrix}.
\end{multline}
\modif{Now write} the problem in terms of $U^*$ with a change of variable, to get
\begin{align}
    d_t U^* + [AU^* + B(U^*) + CU^* + \frac{1}{\rho_0} \mathbf{P} \begin{pmatrix}
        \nabla_H p \\0 \\0
    \end{pmatrix} + F_\sigma(U^*)]dt = G_\sigma(U^*)dW_t - \frac{1}{\rho_0} \mathbf{P} \begin{pmatrix}
        \nabla_H dp_t^\sigma \\0 \\0
    \end{pmatrix},
\end{align}
where the operators $F_\sigma$ and $G_\sigma$ are defined as
\begin{align*}
    F_\sigma(U^*) dt &= \mathbf{P} \Big[ d_t \US + [ B(U^*,\US) - \frac{1}{2} \nabla \bcdot (a \nabla \US)]dt + A\US + \Gamma \US - \frac{1}{2} \nabla \bcdot (a \nabla U^*) dt \Big],\\
    G_\sigma(U^*) dW_t &= \mathbf{P} \Big[ - (\sigma dW_t \bcdot \nabla) U^* - (\sigma dW_t \bcdot \nabla) \US - A (\sigma^H dW_t) - \Gamma (\sigma^H dW_t) \Big],
\end{align*}
with $(v_s,w_s)^\tr = \uS = \frac{1}{2}(\nabla \bcdot a)$, $\US = (v_s, 0, 0)^\tr$. Remind that $\uS$ is divergence free, so $w_s = w(v_s)$ from the definition of operator $w(v)$ \modif{-- that is equation \eqref{w-definition}.} In addition, we can derive the following relations for the pressure terms, using equations \eqref{primitive-p} and \eqref{primitive-dp},
\begin{align*}
    \frac{1}{\rho_0} \nabla_H p = -g \: &\nabla_H \int_z^0 {(\beta_T T + \beta_S S) dz'} \\
    + &\nabla_H \int_z^0 \Big[K* [\uS \bcdot \nabla (w(v^*) + w_s)] - \frac{1}{2} \nabla \bcdot (a^K \nabla (w(v^*) + w_s))\Big] dz'\\
    - \frac{1}{\rho_0} \nabla_H (dp_t^\sigma) = \nabla_H &\int_z^0 \Big[K* [\sigma dW_t \bcdot \nabla (w(v^*) + w_s)] + A (\sigma^w dW_t) \Big]dz' + \frac{1}{\rho_0} \nabla_H dp_t^{\sigma,s}.
\end{align*}
The quantities $p^s \: dt$ and $dp_t^{\sigma,s}$ are respectively the bounded variation and the martingale contributions to the surface pressure. As they are independent on the z-axis, we have for all $v^\sharp \in V_1$,
$$\Big(\frac{1}{\rho_0} \mathbf{P}^v \nabla_H p^s, v^\sharp \Big)_H = -\Big(\frac{1}{\rho_0} p^s, \nabla_H \bcdot \mathcal{A} v^\sharp \Big)_H =0,$$
using the boundary conditions on $\mathcal{A}v^\sharp$, which are $\mathcal{A}v^\sharp =0$ on $\mathcal{S}_H$ and $\mathcal{A}v^\sharp\bcdot\vec{n} = \frac{\partial}{\partial\vec{n}}\mathcal{A}v^\sharp \times \vec{n}=0$ on $\partial \mathcal{S}_H$. This shows that $\mathbf{P}^v \nabla_H p^s = 0$. Similarly, we get $\mathbf{P}^v \nabla_H dp_t^{\sigma,s} = 0$. Therefore, we get the following relation,
\begin{align*}
    \frac{1}{\rho_0} &\mathbf{P}^v [\nabla_H p] =  \mathbf{P}^v \Bigg[ -g \:  \nabla_H \int_z^0 {(\beta_T T + \beta_S S) dz'} \Bigg] \\
    &+ \mathbf{P}^v \Bigg[\nabla_H \int_z^0 \Big[ K* [\uS \bcdot \nabla (w(v^*) + w_s)] - \frac{1}{2} \nabla \bcdot (a^K \nabla (w(v^*) + w_s))\Big] dz' \Bigg],\\
    - \frac{1}{\rho_0} &\mathbf{P}^v[\nabla_H (dp_t^\sigma) ] = \mathbf{P}^v \Bigg[ \nabla_H \int_z^0 \Big[K*[\sigma dW_t \bcdot \nabla (w(v^*) + w_s)] + A (\sigma^w dW_t) \Big]dz'\Bigg].
\end{align*}

\modif{On the one hand}, the bounded variation surface pressure $p^s$ corresponds to a Lagrange multiplier associated with the constraint $\nabla_H \bcdot \mathcal{A} v=0$. This reminds of Cao and Titi's proof in the deterministic setting \cite{CT_2007}, where they showed that, up to some coupling terms, the barotropic mode follows a \modif{2D Navier-Stokes dynamics} while the baroclinic mode follows a 3D Burgers \modif{dynamics.} Crucially, the bounded variation surface pressure does not influence the baroclinic dynamics, given the divergence-free nature of the barotropic mode under vertical boundary conditions. On the other hand, the martingale surface pressure $p^{s,\sigma}$ emerges from our proposed stochastic modelling, acting as a perturbation to the pressure $p^s$. This additional term $p^{s,\sigma}$ \modif{may} impact both the barotropic and baroclinic dynamics. However, under the strong hydrostatic hypothesis, the martingale pressure equation simplifies to $\partial_z dp_t^\sigma = 0$ \modif{from equation} \eqref{strong-hydro}, meaning $dp_t^\sigma = dp_t^{s,\sigma}$. Hence, the martingale pressure term solely affects the barotropic dynamics in this scenario, akin to the deterministic case, enabling the use of analogous methodologies. Our derivations lead to the following formulation of the filtered problem.

\paragraph{Low-pass filtered problem $(\mathcal{P}_K)$:}
For any $K \in H^3(\mathcal{S}, \R)$, we define $(\mathcal{P}_K)$, the abstract primitive equations problem with weak low-pass filtered hydrostatic hypothesis, as follows,
$$\left\{\begin{array}{l}
     d_t U^* + [AU^* + B(U^*) + CU^* + \frac{1}{\rho_0} \mathbf{P} \nabla_H p  + F_\sigma(U^*)]dt = G_\sigma(U^*)dW_t - \frac{1}{\rho_0} \mathbf{P} \nabla_H dp_t^\sigma,\\
     \begin{array}{rcl}
        - \frac{1}{\rho_0} \mathbf{P} \nabla_H (dp_t^\sigma) & = & \mathbf{P} \Bigg[ \nabla_H \int_z^0 \Big[K * [\sigma dW_t \bcdot \nabla (w(v^*) + w_s)] + A (\sigma^w dW_t) \Big]dz' \Bigg],\\
        \frac{1}{\rho_0} \mathbf{P} \nabla_H p & = & \mathbf{P} \Bigg[ \nabla_H \int_z^0 K*\Big[\uS \bcdot \nabla (w(v^*) + w_s) \Big] dz'\\
        & & - \nabla_H \int_z^0 \frac{1}{2} \nabla \bcdot (a^K \nabla (w(v^*) + w_s)) dz' - g \:  \nabla_H \int_z^0 {(\beta_T T + \beta_S S) dz'} \Bigg], 
    \end{array}
\end{array}\right.$$
under the condition $\int_{-h}^0 \nabla_H \bcdot v^* =0$. The problem is supplemented with the boundary conditions \eqref{boundary-conditions} and the initial conditions \eqref{initial-conditions}. As mentioned before, we assume that the noise, its variance tensor and their gradients cancel on the horizontal boundary $\Gamma_l$, and that the vertical noise $\sigma_z dW_t$ cancels on the vertical upper and bottom boundaries $\Gamma_u \cup \Gamma_b$. Moreover, the noise is assumed to follow the regularity conditions \eqref{smoothness-noise} and \eqref{smoothness-noise-2}.

\section{Main results}

Our main results concern the well-posedness of the weak low-pass filtered problem $(\mathcal{P}_K)$,
\begin{theorem} \label{theorem-weak}
Suppose $K \in H^{3}(\mathcal{S}, \R)$. Then, the following propositions hold,
\begin{enumerate}
    \item The problem $(\mathcal{P}_K)$ admits at least one global-in-time martingale solution, for all $T>0$, in the space $$L^2\Big(\Omega,L^2\big([0,T],V\big)\Big) \cap L^2\Big(\Omega,L^\infty\big([0,T],H\big)\Big).$$
    \item There exists a stopping time $\tau>0$ such that $(\mathcal{P}_K)$ admits a local-in-time pathwise solution $U^*$, which fulfils, for all $T>0$ and for all stopping time $0< \tau' < \tau$,
    $$U_{\tau' \wedge \bcdot}^* \in L^2\Big(\Omega, L^2\big([0,T],\mathcal{D}(A)\big)\Big) \cap L^2\Big(\Omega, C([0,T],V)\Big).$$
    This solution is unique up to indistinguishability, that is for all solutions $U^*$ and $\hat{U}^*$ of $(\mathcal{P}_K)$ associated to the stopping times $\tau, \tilde{\tau}$ respectively, the following holds,
    $$\mathbbm{P}\Big(\sup_{[0,T]} \|U_{\tau' \wedge \bcdot}^* - \hat{U}_{\tau' \wedge \bcdot}^*\|_H^2 = 0 ;\quad \forall T>0 \Big)=1,$$
    for all stopping time $0< \tau' < \tau \wedge \hat{\tau}$.
\end{enumerate}
\end{theorem}
In addition, we propose the following assumption on the noise structure,\\
\textbf{Barotropic horizontal noise assumption (BHN)}
\begin{center} \emph{
    $\sigma_H dW_t$ is independent of the variable $z$,\\ i.e. the horizontal  noise is constant over the $z$ axis.}
\end{center}
Such assumption is used to demonstrate our results of global pathwise well-posedness, and continuity with respect to the initial data and the noise data. Namely,
\begin{theorem}
    Assume that $(BHN)$ holds, and choose $K \in H^{15/4}(\mathcal{S}, \R)$. Then, the problem $(\mathcal{P}_K)$ admits a global-in-time pathwise solution, which is unique up to indistinguishability, in the space $$L_{loc}^2\Big([0,+\infty),\mathcal{D}(A)\Big) \cap C\Big([0,+\infty),V\Big).$$ This solution is continuous in the following sense: fix $T>0$, and define $\Sigma$, a space of noise operators, as follows,
    $$\Sigma = \Bigg\{ \sigma \in \mathcal{L}_2\Big(\mathcal{W}, H^4(\mathcal{S}, \R^3)\Big) \: \Bigg| \: a \in H^1\Big([0,T],H^3(\mathcal{S}, \R^{3 \times 3})\Big) \Bigg\}.$$
    Let $(U_0^n, \sigma^n) \in V \times \Sigma$ a sequence converging to $(U_0, \sigma) \in V \times \Sigma$, and denote by $U^n$ and $U$ the solutions associated to $(U_0^n, \sigma^n)$ and $(U_0, \sigma)$, respectively. Then, $U^n \rightarrow U$ in probability, in the space $L^2\big([0,T],\mathcal{D}(A)\big) \cap C([0,T],V)$.   
    In particular, for a fixed initial data, if we denote by $U^{\sigma}$ the solution of $(\mathcal{P}_K)$ associated to the noise data $\sigma \in \mathcal{N}$, then, in probability,
    \begin{equation*}
        U^{\Upsilon^{1/2} \sigma} \underset{ \Upsilon \rightarrow 0}{\longrightarrow} U^0.
    \end{equation*}
    Here, $U^0$ denotes the solution to the problem $(\mathcal{P}_K)$ with noise zero -- i.e. $\sigma = 0$ -- that is to say the classical deterministic primitive equations.
\end{theorem}

We only give a sketch of the proof for Theorem 2 below. Theorem 1 can be obtained by adapting the work proposed in \cite{DHM_2023} to the primitive equations, which can be considered as a simplification of the 3D Navier-Stokes equations. It relies on considering a Galerkin approximation of the problem, then using energy estimates and tightness arguments to show that these solutions converge toward a solution of the initial problem, see \cite{BS_2021, DGHT_2011, DHM_2023, FG_1995}.

Fix $T>0$. Remind that the barotropic and baroclinic modes of velocity are defined in \modif{equation} \eqref{def-barotropic-baroclinic}. Assuming that Theorem 1 hold, we show that the following estimates hold for any stopping times $0 < \eta < \zeta < T$, there exist three constants $C_1, C_2, C_3 > 0$ such that:
\begin{itemize}
    \item Barotropic velocity estimate in ``$H^1-H^2$''
        \begin{multline}
        \E\Big[\|\bar{v}\|_{V}^2(\zeta)\Big] + \E\Big[ \int_\eta^\zeta \|A \bar{v}\|_H^2 \: ds\Big]
        \leq C_1 \E\Big[\|\bar{v}(\eta)\|_{V}^2 + 1\Big] \label{ito-barotrope-simple}\\
        + C_1 \E\Big[\int_\eta^\zeta \|U\|_H^2 \|U\|_V^2 \|\bar{v}\|_{V}^2 + \|U\|_V^2 + \int_\mathcal{S} |\tilde{v}|^2 |\nabla_3 \tilde{v}|^2 + \|\partial_z v\|_V^2 ds\Big],
        \end{multline}
    \item Vertical gradient of velocity estimate in ``$L^2-H^1$''
        \begin{multline}
            \E\Big[\sup_{[\eta,\zeta]} \|\partial_z v\|_{H}^2\Big] + \E\Big[ \int_\eta^\zeta \|\partial_z v\|_V^{2} \: ds\Big]
            \leq C_2 \E\Big[\|\partial_z v(\eta)\|_{H}^2 +1\Big] \label{estimation-gradient-z} \\
            + C_2 \E\Big[\int_\eta^\zeta (1+\|U\|_V^2) (1+\|\partial_z v\|_H^2) + \int_\mathcal{S} |\nabla \tilde{v}|^2 |\tilde{v}|^2 ds \Big].
        \end{multline}
    \item Baroclinic velocity estimate in ``$L^4$''
        \begin{multline}
            \E\Big[\sup_{[\eta,\zeta]} |\tilde{v}\|_{L^4}^4\Big] + \E\Big[ \int_\eta^\zeta \int_{\mathcal{S}} (|\nabla \tilde{v}|^2 |\tilde{v}|^2 + |\nabla |\tilde{v}|^2|^2 + |\partial_z \tilde{v}|^2 |\tilde{v}|^2 + |\partial_z |\tilde{v}|^2 |^2) d\mathcal{S} \: ds\Big] \label{estimation-barocline}\\
            \leq C_3 \E\Big[1+ \|\tilde{v}(\eta)\|_{L^4}^4 + \int_\eta^\zeta (\|U\|_V^2 +1)\|\tilde{v}\|_{L^4}^4  \Big] + \frac{1}{16(C_2 \vee 1) } \E\Big[\int_\eta^\zeta \|\partial_z v\|_V^2 \: ds \Big].
        \end{multline}
\end{itemize}
Once gathered, these estimates lead to the existence of a global pathwise solution to the problem $(\mathcal{P}_K)$, using similar arguments as in \cite{AHHS_2022} \modif{-- in particular the stochastic Grönwall's lemma, see \cite{GHZ_2009}.} Moreover, using similar arguments, we may show the uniqueness of the solution, as well as the continuity of the solution with respect to the initial condition and the noise operator $\sigma$. These estimates are related to the proof of the well-posedness of the deterministic primitive equations proposed originally in the article \cite{CT_2007}. Essentially, the argument is that $\bar{v}$ follows a 2D Navier Stokes equation, while $\tilde{v}$ follows a 3D Burger equation, up to some coupling terms. Therefore, the idea is to estimate $\bar{v}$ in the strong sense (``$H^1 - H^2$'') and $\tilde{v}$ in an $L^p$-space -- namely $L^4$. Moreover, a third estimate is needed on the vertical velocity gradient $\partial_z v$ (``$L^2 - H^1$''). In addition, the noise structure constraints $(BHN)$ we applied for ensuring global existence remind of those \modif{proposed} in \cite{AHHS_2022, AHHS_2023_preprint}, \modif{particularly} $\sigma_H dW_t = (\sigma_H dW_t) (x,y)$ and $\sigma_z dW_t = (\sigma_z dW_t) (x,y,z)$.

However, the scheme of proof in \cite{AHHS_2022} differs from ours, which is closer in spirit to the one of \cite{BS_2021}. Notably, our proof accommodates scenarios where the vertical acceleration $\mathbbm{D}_t w$ is entirely neglected, which corresponds to choosing $K=0$ and neglecting additive noise in \modif{equation} \eqref{weak-split}. Our well-posedness outcome is similar as the one of \cite{AHHS_2022}, except that we transitioned from water world to rigid-lid boundary conditions. Nonetheless, assuming that the noise and its gradient cancel on the boundary remains pivotal for the validity of our integration by parts arguments.

\modif{A more detailed proof will be submitted subsequently.}

\paragraph{Remark:} The barotropic horizontal noise assumption $(BHN)$ aligns with the validity domain of the primitive equations in their deterministic form. These equations hold true when the squared aspect ratio $\alpha^2 := (h/L)^2$ is negligible compared to the Richardson number $Ri := \frac{N^2}{(\partial_z v)^2}$ \cite{MHPA_1997}. Here, $v$ denotes horizontal velocity, $h$ the ocean depth, $L$ the horizontal scale (e.g., $\sqrt{|\mathcal{S}_H|}$), and $N^2= - \frac{g}{\rho_0} \partial_z \rho$ the Brünt-Väisälä frequency. This condition reads
\begin{align}
    \frac{\alpha^2}{Ri} \ll 1 \quad \text{or equivalently} \quad (\partial_z v)^2 \ll \frac{N^2}{\alpha^2}. \label{condition-Marshall}
\end{align}
This particularly holds in the limit of small vertical shear of the horizontal component. In such case, the horizontal velocity becomes almost independent of $z$, which we call ``quasi-barotropic''. In stochastic flow contexts, the horizontal noise models a small-scale velocity, denoted by $\Upsilon^{1/2} v'$, where $\Upsilon^{1/2}$ is a scaling factor ensuring $v$ and $v'$ share the same order of magnitude. Hence, condition \eqref{condition-Marshall} can be rewritten as
\begin{align}
(\partial_z (v + \Upsilon^{1/2} v'))^2 \ll \frac{N^2}{\alpha^2}, \quad \text{or} \quad (\partial_z v)^2 , \: \Upsilon (\partial_z v')^2 \ll \frac{N^2}{\alpha^2}. \label{condition-Marshall-stocha}
\end{align}
Consequently, the LU stochastic primitive equations remain physically valid under condition \eqref{condition-Marshall-stocha}, \modif{that is} when the horizontal noise modeling $\Upsilon^{1/2} v'$ is either sufficiently small ($\Upsilon \rightarrow 0$) or quasi-barotropic ($(\partial_z v')^2 \rightarrow 0$). In this setup, the noise structure hypothesis in Theorem 1.3 is equivalent to $\partial_z v' = 0$. We anticipate that a slight deviation from this assumption -- \modif{that is to say} considering a noise with a non-zero small \modif{enough} baroclinic mode -- would yield similar well-posedness results. This is because the energy stemming from the baroclinic noise -- when small enough -- is likely to be balanced by the combined molecular and stochastic diffusions. However, considering a large baroclinic noise component seems to be a serious challenge when proving the LU primitive equations well-posedness.

\printbibliography[title=References]

@article{DHM_2023,
	author = {Debussche, Arnaud and Hug, B{\'e}renger and M{\'e}min, Etienne},
	doi = {10.1007/s00021-023-00764-0},
	id = {Debussche2023},
	isbn = {1422-6952},
	journal = {Journal of Mathematical Fluid Mechanics},
	number = {1},
	pages = {19},
	title = {A Consistent Stochastic Large-Scale Representation of the Navier--Stokes Equations},
	url = {https://doi.org/10.1007/s00021-023-00764-0},
	volume = {25},
	year = {2023},
    month ={1}}

@article{DGHT_2011,
	doi = {10.1016/j.physd.2011.03.009},
  
	url = {https://doi.org/10.1016%2Fj.physd.2011.03.009},
  
	year = 2011,
	month = {7},
  
	publisher = {Elsevier {BV}},
  
	volume = {240},
  
	number = {14-15},
  
	pages = {1123--1144},
  
	author = {Arnaud Debussche and Nathan Glatt-Holtz and Roger Temam},
  
	title = {Local martingale and pathwise solutions for an abstract fluids model},
  
	journal = {Physica D: Nonlinear Phenomena}
}

@article{BS_2021,
	author = {Zdzislaw Brze{{\'z}}niak and Jakub Slav{\'i}k},
	doi = {https://doi.org/10.1016/j.jde.2021.05.049},
	issn = {0022-0396},
	journal = {Journal of Differential Equations},
	pages = {617-676},
	title = {Well-posedness of the 3D stochastic primitive equations with multiplicative and transport noise},
	url = {https://www.sciencedirect.com/science/article/pii/S0022039621003521},
	volume = {296},
	year = {2021}}

@article{FG_1995,
	title = {Martingale and stationary solutions for stochastic {Navier}-{Stokes} equations},
	volume = {102},
	issn = {1432-2064},
	url = {https://doi.org/10.1007/BF01192467},
	doi = {10.1007/BF01192467},
	number = {3},
	journal = {Probability Theory and Related Fields},
	author = {Flandoli, Franco and Gatarek, Dariusz},
	month = sep,
	year = {1995},
	pages = {367--391}
}

@article{AHHS_2022,
	author = {Antonio Agresti and Matthias Hieber and Amru Hussein and Martin Saal},
	doi = {10.1007/s40072-022-00277-3},
	id = {Agresti2022},
	isbn = {2194-041X},
	journal = {Stochastics and Partial Differential Equations: Analysis and Computations},
	title = {The stochastic primitive equations with transport noise and turbulent pressure},
	url = {https://doi.org/10.1007/s40072-022-00277-3},
	year = {2022},
    month= {10}}

@misc{AHHS_2023_preprint,
      title={The stochastic primitive equations with non-isothermal turbulent pressure}, 
      author={Antonio Agresti and Matthias Hieber and Amru Hussein and Martin Saal},
      year={2023},
      eprint={2210.05973},
      archivePrefix={arXiv},
      primaryClass={math.AP}
}

@article{Mémin_2014,
  title={Fluid flow dynamics under location uncertainty},
  author={M{\'e}min, Etienne},
  journal={Geophysical \& Astrophysical Fluid Dynamics},
  volume={108},
  number={2},
  pages={119--146},
  year={2014},
  publisher={Taylor \& Francis}
}

@inproceedings{TML_2023,
	address = {Cham},
	title = {Primitive {Equations} {Under} {Location} {Uncertainty}: {Analytical} {Description} and {Model} {Development}},
	isbn = {978-3-031-18988-3},
	booktitle = {Stochastic {Transport} in {Upper} {Ocean} {Dynamics}},
	publisher = {Springer International Publishing},
	author = {Tucciarone, Francesco L and Mémin, Etienne and Li, Long},
	editor = {Chapron, Bertrand and Crisan, Dan and Holm, Darryl and Mémin, Etienne and Radomska, Anna},
	year = {2023},
	pages = {287--300}
}

@article{CT_2007,
  title={Global well-posedness of the three-dimensional viscous primitive equations of large scale ocean and atmosphere dynamics},
  author={Cao, Chongsheng and Titi, Edriss S},
  journal={Annals of Mathematics},
  pages={245--267},
  year={2007},
  publisher={JSTOR}
}

@article{Berner-Al_2017,
  title={Stochastic parameterization: Toward a new view of weather and climate models},
  author={Berner, Judith and Achatz, Ulrich and Batte, Lauriane and Bengtsson, Lisa and De La Camara, Alvaro and Christensen, Hannah M and Colangeli, Matteo and Coleman, Danielle RB and Crommelin, Daan and Dolaptchiev, Stamen I and others},
  journal={Bulletin of the American Meteorological Society},
  volume={98},
  number={3},
  pages={565--588},
  year={2017},
  publisher={American Meteorological Society}
}

@article{FOBWL_2015,
  title={Stochastic climate theory and modeling},
  author={Franzke, Christian LE and O'Kane, Terence J and Berner, Judith and Williams, Paul D and Lucarini, Valerio},
  journal={Wiley Interdisciplinary Reviews: Climate Change},
  volume={6},
  number={1},
  pages={63--78},
  year={2015},
  publisher={Wiley Online Library}
}

@book{FO_2017,
  title={Nonlinear and stochastic climate dynamics},
  author={Franzke, Christian LE and O'Kane, Terence J},
  year={2017},
  publisher={Cambridge University Press}
}

@article{MTV_1999,
  title={Models for stochastic climate prediction},
  author={Majda, Andrew J and Timofeyev, Ilya and Vanden Eijnden, Eric},
  journal={Proceedings of the National Academy of Sciences},
  volume={96},
  number={26},
  pages={14687--14691},
  year={1999},
  publisher={National Acad Sciences}
}

@article{BCCLM_2020,
  title={Deciphering the role of small-scale inhomogeneity on geophysical flow structuration: a stochastic approach},
  author={Bauer, Werner and Chandramouli, Pranav and Chapron, Bertrand and Li, Long and M{\'e}min, Etienne},
  journal={Journal of Physical Oceanography},
  volume={50},
  number={4},
  pages={983--1003},
  year={2020},
  publisher={American Meteorological Society}
}

@article{RMC_2017_Pt1,
  title={Geophysical flows under location uncertainty, Part I Random transport and general models},
  author={Resseguier, Valentin and M{\'e}min, Etienne and Chapron, Bertrand},
  journal={Geophysical \& Astrophysical Fluid Dynamics},
  volume={111},
  number={3},
  pages={149--176},
  year={2017},
  publisher={Taylor \& Francis}
}

@article{RMC_2017_Pt2,
  title={Geophysical flows under location uncertainty, part II quasi-geostrophy and efficient ensemble spreading},
  author={Resseguier, Valentin and M{\'e}min, Etienne and Chapron, Bertrand},
  journal={Geophysical \& Astrophysical Fluid Dynamics},
  volume={111},
  number={3},
  pages={177--208},
  year={2017},
  publisher={Taylor \& Francis}
}

@article{RMC_2017_Pt3,
  title={Geophysical flows under location uncertainty, Part III SQG and frontal dynamics under strong turbulence conditions},
  author={Resseguier, Valentin and M{\'e}min, Etienne and Chapron, Bertrand},
  journal={Geophysical \& Astrophysical Fluid Dynamics},
  volume={111},
  number={3},
  pages={209--227},
  year={2017},
  publisher={Taylor \& Francis}
}

@article{RMHC_2017,
  title={Stochastic modelling and diffusion modes for proper orthogonal decomposition models and small-scale flow analysis},
  author={Resseguier, Valentin and M{\'e}min, Etienne and Heitz, Dominique and Chapron, Bertrand},
  journal={Journal of Fluid Mechanics},
  volume={826},
  pages={888--917},
  year={2017},
  publisher={Cambridge University Press}
}

@article{RPMC_2021,
  title={Quantifying truncation-related uncertainties in unsteady fluid dynamics reduced order models},
  author={Resseguier, Valentin and Picard, Agustin M and M{\'e}min, Etienne and Chapron, Bertrand},
  journal={SIAM/ASA Journal on Uncertainty Quantification},
  volume={9},
  number={3},
  pages={1152--1183},
  year={2021},
  publisher={SIAM}
}

@article{TCM_2021,
  title={Stochastic linear modes in a turbulent channel flow},
  author={Tissot, Gilles and Cavalieri, Andr{\'e} VG and M{\'e}min, Etienne},
  journal={Journal of Fluid Mechanics},
  volume={912},
  pages={A51},
  year={2021},
  publisher={Cambridge University Press}
}

@article{CHLM_2018,
  title={Coarse large-eddy simulations in a transitional wake flow with flow models under location uncertainty},
  author={Chandramouli, Pranav and Heitz, Dominique and Laizet, Sylvain and M{\'e}min, Etienne},
  journal={Computers \& Fluids},
  volume={168},
  pages={170--189},
  year={2018},
  publisher={Elsevier}
}

@article{CMH_2020,
  title={4D large scale variational data assimilation of a turbulent flow with a dynamics error model},
  author={Chandramouli, Pranav and M{\'e}min, Etienne and Heitz, Dominique},
  journal={Journal of Computational Physics},
  volume={412},
  pages={109446},
  year={2020},
  publisher={Elsevier}
}

@article{HM_2017,
  title={Stochastic representation of the Reynolds transport theorem: revisiting large-scale modeling},
  author={Harouna, Souleymane K and M{\'e}min, Etienne},
  journal={Computers \& Fluids},
  volume={156},
  pages={456--469},
  year={2017},
  publisher={Elsevier}
}

@article{CL_1976,
  title={A rational model for Langmuir circulations},
  author={Craik, Alex DD and Leibovich, Sidney},
  journal={Journal of Fluid Mechanics},
  volume={73},
  number={3},
  pages={401--426},
  year={1976},
  publisher={Cambridge University Press}
}

@article{GOP_2004,
  title={Mathematical perspectives on large eddy simulation models for turbulent flows},
  author={Guermond, Jean-Luc and Oden, John T and Prudhomme, Serge},
  journal={Journal of Mathematical Fluid Mechanics},
  volume={6},
  pages={194--248},
  year={2004},
  publisher={Springer}
}

@article{MSM_1997,
  title={Langmuir turbulence in the ocean},
  author={McWilliams, James C and Sullivan, Peter P and Moeng, Chin-Hoh},
  journal={Journal of Fluid Mechanics},
  volume={334},
  pages={1--30},
  year={1997},
  publisher={Cambridge University Press}
}

@article{GHZ_2009,
author = {Nathan Glatt-Holtz and Mohammed Ziane},
title = {{Strong pathwise solutions of the stochastic Navier-Stokes system}},
volume = {14},
journal = {Advances in Differential Equations},
number = {5/6},
publisher = {Khayyam Publishing, Inc.},
pages = {567 -- 600},
year = {2009},
doi = {10.57262/ade/1355867260},
URL = {https://doi.org/10.57262/ade/1355867260}
}

@book{Vallis2017,
  title={Atmospheric and oceanic fluid dynamics},
  author={Vallis, Geoffrey K},
  year={2017},
  publisher={Cambridge University Press}
}

@article{CFH_2019,
	author = {Crisan, Dan and Flandoli, Franco and Holm, Darryl},
	journal = {Journal of Nonlinear Science},
	number = {3},
	pages = {813--870},
	publisher = {Springer},
	title = {Solution properties of a 3{D} stochastic {E}uler fluid equation},
	volume = {29},
	year = {2019}}

@article{GCL_2023,
  title={Existence and uniqueness of maximal solutions to SPDEs with applications to viscous fluid equations},
  author={Goodair, Daniel and Crisan, Dan and Lang, Oana},
  journal={Stochastics and Partial Differential Equations: Analysis and Computations},
  pages={1--64},
  year={2023},
  publisher={Springer}
}

@article{LCM_2023,
	author = {Lang, Oana and Crisan, Dan and M{\'e}min, Etienne},
	date = {2023-02-20},
	doi = {10.1007/s00021-023-00769-9},
	id = {Lang2023},
	isbn = {1422-6952},
	journal = {Journal of Mathematical Fluid Mechanics},
	number = {2},
	pages = {29},
	title = {Analytical Properties for a Stochastic Rotating Shallow Water Model Under Location Uncertainty},
	url = {https://doi.org/10.1007/s00021-023-00769-9},
	volume = {25},
	year = {2023}}

@article{MR_2005,
	author = {Remigijus Mikulevicius and Boris L Rozovskii},
	date = {2005-01-01},
	doi = {10.1214/009117904000000630},
	journal = {The Annals of Probability},
	journal1 = {The Annals of Probability},
	journal2 = {The Annals of Probability},
	month = {1},
	number = {1},
	pages = {137--176},
	title = {Global {L}$_{2}$-solutions of stochastic {N}avier--{S}tokes equations},
	url = {https://doi.org/10.1214/009117904000000630},
	volume = {33},
	year = {2005}}

@article{FGP_10,
	author = {Flandoli, Franco and Gubinelli, Massimiliano and Priola, Enrico },
	date = {2010-04-01},
	doi = {10.1007/s00222-009-0224-4},
	id = {Flandoli2010},
	isbn = {1432-1297},
	journal = {Inventiones mathematicae},
	number = {1},
	pages = {1--53},
	title = {Well-posedness of the transport equation by stochastic perturbation},
	url = {https://doi.org/10.1007/s00222-009-0224-4},
	volume = {180},
	year = {2010}}

@article{FL_2021,
	author = {Flandoli, Franco and Luo, Dejun},
	date = {2021-06-01},
	doi = {10.1007/s00440-021-01037-5},
	id = {Flandoli2021},
	isbn = {1432-2064},
	journal = {Probability Theory and Related Fields},
	number = {1},
	pages = {309--363},
	title = {High mode transport noise improves vorticity blow-up control in 3D Navier--Stokes equations},
	url = {https://doi.org/10.1007/s00440-021-01037-5},
	volume = {180},
	year = {2021}}

@article{FGL_2021,
	author = {Franco Flandoli and Lucio Galeati and Dejun Luo},
	doi = {10.1080/03605302.2021.1893748},
	eprint = {https://doi.org/10.1080/03605302.2021.1893748},
	journal = {Communications in Partial Differential Equations},
	number = {9},
	pages = {1757-1788},
	publisher = {Taylor & Francis},
	title = {Delayed blow-up by transport noise},
	url = {https://doi.org/10.1080/03605302.2021.1893748},
	volume = {46},
	year = {2021}}

@article{KT_2023,
      title={Global Well-Posedness of the Primitive Equations of Large-Scale Ocean Dynamics with the Gent-McWilliams-Redi Eddy Parametrization Model}, 
      author={Peter Korn and Edriss S Titi},
      year={2023},
      eprint={2304.03242},
      archivePrefix={arXiv},
      primaryClass={math.AP}
}

@article{MHPA_1997,
author = {Marshall, John and Hill, Chris and Perelman, Lev and Adcroft, Alistair},
title = {Hydrostatic, quasi-hydrostatic, and nonhydrostatic ocean modeling},
journal = {Journal of Geophysical Research: Oceans},
volume = {102},
number = {C3},
pages = {5733-5752},
doi = {https://doi.org/10.1029/96JC02776},
url = {https://agupubs.onlinelibrary.wiley.com/doi/abs/10.1029/96JC02776},
eprint = {https://agupubs.onlinelibrary.wiley.com/doi/pdf/10.1029/96JC02776},
year = {1997}
}

@book{book_DPZ_2014,
  title={Stochastic equations in infinite dimensions},
  author={Da Prato, Giuseppe and Zabczyk, Jerzy},
  year={2014},
  publisher={Cambridge university press, Second Edition}
}

@book{book_GM_2010,
  title={Stochastic differential equations in infinite dimensions: with applications to stochastic partial differential equations},
  author={Gawarecki, Leszek and Mandrekar, Vidyadhar},
  year={2010},
  publisher={Springer Science \& Business Media}
}

\end{document}